\input amstex 
\documentstyle{amsppt} 
\NoBlackBoxes 
\NoRunningHeads 
\magnification=\magstep1 
\def\cit#1#2{\ifx#1!\cite{#2}\else#2\fi} 
\define\g{\frak g} 
\topmatter 
\title 
Description of infinite dimensional abelian regular Lie groups  
\endtitle 
\author 
Peter Michor \\ 
Josef Teichmann 
\endauthor 
\address 
J\. Teichmann: 
Institut f\"ur Mathematik, Universit\"at Wien 
Strudlhofgasse 4, A-1090 Wien, Austria 
\endaddress 
\email 
jteichma\@radon.mat.univie.ac.at 
\endemail 
\address 
P\. Michor: 
Institut f\"ur Mathematik, Universit\"at Wien 
Strudlhofgasse 4, A-1090 Wien, Austria; and: 
Erwin Schr\"odinger Institut f\"ur Mathematische Physik, 
Boltzmanngasse, A-1090 Wien, Austria 
\endaddress 
\email peter.michor\@esi.ac.at \endemail 
\thanks  Supported  by `Fonds zur F\"orderung der wissenschaftlichen   
Forschung, Projekt P~10037~PHY'. \endthanks \abstract 
It is shown that every abelian regular Lie group is a quotient of its  
Lie algebra via the exponential mapping. 
\endabstract 
\endtopmatter 
\document 
This paper is a sequel of \cit!{3}, see also \cit!{4}, chapter VIII,  
where a regular Lie group is defined as a smooth Lie group modeled on  
convenient vector spaces such that the right logarithmic derivative  
has a smooth inverse  
$\operatorname{Evol}: C^\infty(\Bbb R,\frak g) \to C^\infty(\Bbb R,G)$, 
the canonical evolution operator, where $\frak g$ is the Lie algebra.  
We follow the notation and the concepts of this paper closely. 
 
\proclaim{Lemma} Let $G$ be an abelian regular Lie group with Lie  
algebra $\frak \g$. Then the evolution operator is given by  
$ 
\text{\rm Evol}(X)(t) := {\text{\rm Evol}}^r(X) (t) =  
\exp \bigl( \int_0^t X(s) ds \bigr)  
$ 
for $ X \in C^{\infty} ({\Bbb R},\frak{g}) $. 
\endproclaim 
 
\demo{Proof} Since $G$ is regular it has an exponential mapping 
$ \exp: \frak{g} \to G $ which is a smooth group 
homomorphism, because $s\mapsto \exp(sX)\exp(sY) $ is a 
smooth one-parameter group in $G$ with generator $X + Y$, thus  
$\exp(X) \exp(Y) = \exp(X + Y)$ by uniqueness, \cit!{3}, 3.6 or  
\cit!{4}, 36.7. 
The Lie algebra $\g$ is a convenient vector space with  
evolution mapping  
$ 
\operatorname{Evol}_{\g}(X)(t) = \int_0^1 X(s)ds,  
$ 
see \cit!{3}, 5.4, or \cit!{4}, 38.5. The mapping $\exp:\g\to G$ is  
a homomorphism of Lie groups and thus intertwines the evolution  
operators by \cit!{3}, 5.3 or \cit!{4}, 38.4, hence the formula. 
 
Another proof is by differentiating the right hand side, using  
\cit!{3}, 5.10 or \cit!{4}, 38.2. 
\qed\enddemo 
 
As consequence we obtain that an abelian Lie group $G$ is regular if  
and only if an exponential map exists. 
Furthermore, an exponential map is surjective on a 
connected abelian Lie 
group, because $ \exp(\int_0^t \delta^r c(s) ds) =  
\text{\rm Evol} ({\delta}^r c)(t)  = c(t)$  
for any smooth curve $c : {\Bbb R} \to G$ with $c(0) = e$. 
 
\proclaim{Theorem} Let $ G $ be an abelian, connected 
and regular Lie group, then there is a $c^{\infty} $-open neighborhood  
$V$ of zero in $ \frak{g} $ so that $ \exp(V) $ is open in $G$ and  
$ \exp: V \to \exp(V) $ is a diffeomorphism.   
Moreover, $\frak{g} / \ker(\exp)\to G$ is an isomorphism of Lie 
groups.  
\endproclaim 
 
\demo{Proof} Given a connected, abelian and regular Lie group $ G $,  
we look at the 
universal covering group $ \tilde{G} @>\pi>> G $, see \cit!{4},  
27.14, which is also abelian and regular. Any tangent Lie algebra  
homomorphism from a simply connected Lie group to a regular Lie group  
can be uniquely integrated to a Lie group homomorphism by \cit!{5} or  
\cit!{3}, 7.3 or \cit!{4}, 40.3. Consequently, there exists a  
homomorphism $\Phi:\tilde{G}\to\frak{g}$ with $\Phi'=id_{\frak{g}}$.  
Since $\tilde{G}$ is regular there is a map from $\frak{g}$ to  
$\tilde{G}$ extending $id$, which has to be the inverse of $\phi$  
and which is a fortiori the exponential map $\widetilde{\exp}$ of  
$\tilde{G}$, so $\Phi$ is an isomorphism of Lie  
groups. The universal covering projection  $\pi$ intertwines  
$\widetilde\exp$ and $\exp$, so the result follows. The quotient  
$\frak{g}/\ker(\exp)$ is a Lie group since there are natural chart  
maps and the quotient space is a Hausdorff space by the Hausdorff  
property on $G$. \qed 
\enddemo 
 
\subhead Remarks \endsubhead 
Given a convenient vector space $E$ and a subgroup $Z$, it is 
not obvious how to determine simple conditions to ensure that $ E/Z $ is a 
Hausdorff space, because $ c^{\infty}E $ is not a topological vector space 
in general(see \cit!{4}, Chapter I):  
An additive subgroup $Z$ of $ E $ is called 'discrete' if  
there is a $c^{\infty}$-open zero neighborhood $V$ with  
$V \cap (Z + V)= \{0\}$ and  
for any $x \notin Z$ there is a $c^{\infty}$-open zero neighborhood $U$  
so that $ (x + Z + U) \cap (Z + U) = \emptyset$. The above kernel of  
$\exp$ naturally has this property, consequently any regular connected 
abelian Lie group is a convenient vector space modulo a 'discrete' subgroup. 
 
Let $E$ be a Fr\'echet space, then a subgroup is 'discrete' if and only if 
there is an open zero neighborhood $V$ with $V \cap (Z + V)= \{0\}$,  
because $c^{\infty}E = E$. This leads immediately to a 
generalization of a result of Galanis (\cit!{2}), 
who proved that every abelian Fr\'echet-Lie group which admits 
an exponential map being a local diffeomorphism around zero is a projective 
limit of Banach Lie groups. With the above theorem one can easily write down 
this limit in general. 

With the above methods it is necessary to assume regularity: 
Otherwise one obtains as image of $\Phi$ a dense arcwise connected 
subgroup of the convenient vector space $\g$, which does not allow 
any conclusion in contradiction the finite dimensional case.   
Note that the closed subgroup of  
integer-valued functions in $ L^2([0,1],{\Bbb R}) $  
is arcwise connected but not a Lie subgroup 
(see 
\cit!{1}) so that Yamabe's   
theorem is already 
wrong on the level of infinite dimensional Hilbert spaces.

\enddocument